\input amstex
\documentstyle{amsppt}
\magnification 1200
\TagsOnRight
\NoBlackBoxes
\nologo  
\vsize=8.8 truein
\NoBlackBoxes
\define\refBonaP    {1}
\define\refBonaS    {2}
\define\refCanosa   {3}
\define\refChen     {4}
\define\refCoquel   {5}
\define\refCraig    {6}
\define\refDiPerna  {7}
\define\refGoodman  {8}
\define\refJacobs   {9}
\define\refJeffrey  {10}
\define\refJohnson  {11}
\define\refKato     {12}
\define\refKenigOne {13}
\define\refKenigTwo {14}
\define\refKenigThree  {15}
\define\refKruzkov     {16}
\define\refKruzkovF    {17}
\define\refLaxOne      {18}
\define\refLaxTwo      {19}
\define\refLaxL     {20}
\define\refLaxW     {21}
\define\refMarcati  {22}
\define\refOleinik  {23}
\define\refVonNeumann  {24}
\define\refRaviart     {25}
\define\refRichtmyer   {26}
\define\refSchonbek    {27}
\define\refSzepessy    {28}
\define\refTartar      {29}
\define\refVenakidesOne{30}
\define\refVenakidesTwo{31}

\newif\iftitle

\def\leaderfill{ \leaders\hbox to1em{\hss.\hss}\hfill}

\newwrite\ind
\immediate\openout\ind=indice
\immediate\write\ind{\string\input }
\immediate\write\ind{fmtind}

\def\beginsection#1#2#3\par{\vskip0pt plus.3\vsize\penalty10000
    \vskip0pt plus -.3\vsize
\ifcase#1
\bigskip\vskip\parskip 
\item{\bf #2}  \bf #3\smallskip\or
\noindent{\item{\bf#2}\it #3\hfil\break\vskip0pt}
\else 
{\item{\it#2} \it #3}\fi
\message{#2#3}
\edef\mywrite##1##2##3{%
  \noexpand\writeind {##1}{##2}{##3} }%
\mywrite{#1}{#2}{#3}%
\ifcase#1\mark{#2\quad#3}\fi
\rm}

\def\writeind #1#2#3 {
\write\ind {
 \string\llap{\hbox to30pt{#2\string\hfil}}
                 #3\string\leaderfill \string\quad
 {\string\oldnos \folio}\string\break }
 }

\def\label#1{\smallskip\par\noindent\hangindent30pt #1}

\def\t{\text}
\def\R{\text{{\rm I}\!{\rm R}}} 
\def\N{\Bbb N}

\def\pa{\partial}
\def\ra{\rangle}
\def\la{\langle}

\def\var{\varepsilon}
\def\d{\,\text{\rm  d}}

\def\loc{\,\text{\rm loc}\,}
\def\sgn{\,\text{\rm sgn}\,}
\def\supp{\,\text{\rm supp}\,}

\def\label#1{\smallskip\par\noindent\hangindent30pt#1}

\topmatter 
\title Conservation Laws with Vanishing \\
         Nonlinear Diffusion and Dispersion\endtitle
\thanks Published in : Nonlinear Analysis 36 (1999), 213--230.
\newline This work was partially carried out in the Summer 1992 
during a visit of the first
author at the Istituto per le Applicazioni del Calcolo\endthanks
\endtopmatter
\bigskip
\centerline{\smc Philippe G. LeFloch\footnote"$^1$"{E-mail address: 
LeFloch\@ann.jussieu.fr.}}\medskip

{\baselineskip12pt 

\centerline{{\bf Current address} : Laboratoire Jacques-Louis Lions}
\centerline{Centre National de la Recherche Scientifique}
\centerline{Universit\'e de Paris 6}
\centerline{4, Place Jussieu, 75252 Paris, France.}
 \medskip
\centerline{and}
\bigskip
\centerline{\smc Roberto Natalini\footnote"$^2$"{E-mail address: 
Natalini\@vaxiac.iac.rm.cnr.it}
}}\medskip

{\baselineskip12pt 
\centerline{Istituto per le Applicazioni del Calcolo ``M. Picone''}
\centerline{Consiglio Nazionale delle Ricerche}
\centerline{Viale del Policlinico  137}
\centerline{00161 Roma, Italy}}\medskip

\bigskip\bigskip{\narrower\smallskip\noindent
\baselineskip12pt\eightpoint{\it ABSTRACT\/}. We
study the limiting behavior of the solutions to a class of conservation laws
with vanishing nonlinear diffusion and dispersion terms. We prove the
convergence to the entropy solution of the first order problem 
under a condition on 
the relative size of the diffusion and the dispersion terms. This work is
motivated by the pseudo-viscosity approximation introduced 
by Von Neumann in the 50's.\medskip

\noindent{\it Key words and phrases:}  Nonlinear dispersive waves, 
Korteweg-de Vries equation, pseudo-viscosity,
Burgers equation, shock waves, measure-valued solutions.
\smallskip}

\vfill
\eject 

\beginsection0{1.}{Introduction}

This paper is concerned with the convergence of smooth solutions 
$u=u^{\var,\delta}$
to the initial value problem for the nonlinear dispersive equation
$$
\pa_tu+\pa_x\,f(u)=\var \, \pa_x\beta(\pa_x u)-\delta \, \pa^3_xu,\quad
(x,t)\in{\R}\times(0,\infty),
\tag1.1
$$
with the initial condition
$$
u(x,0)=u^{\var,\delta}_0(x),\quad x\in{\R},
\tag1.2
$$
where the parameters $\var>0$ and $\delta>0$ will converge to zero
and $u^{\var,\delta}_0$ is an approximation of a given initial condition $u_0: 
{\R} \to {\R}$. 
The flux $f=f(u)$ and the (degenerate) viscosity $\beta=\beta(\lambda)$ are
given smooth functions satisfying certain assumptions to be listed shortly. 
When $\var=0$, the equation (1.1) reduces to the generalized Korteweg-de Vries 
(KdV) equation, the original one corresponding to the special flux 
$f(u)= u^2/2$. 
On the other hand, when  $\delta=0$, 
we recover a nonlinear degenerate parabolic equation. 
The first case was extensively studied from mathematical 
and numerical standpoints; 
see for instance \cite{\refKenigOne, \refKenigTwo, \refCraig}
and the references therein (see also Section 2 for some background).
For $\delta=0$ and under suitable assumptions,  
this equation was treated
in \cite{\refMarcati} as a simplified model of the pseudo viscosity 
approximation proposed 
by von Neumann and Richtmyer \cite{\refVonNeumann}
and then studied in \cite{\refRichtmyer, \refRaviart} for numerical purposes.

In this paper we prove, under an assumption on the relative size of the
parameters $\var$ and $\delta$, that the sequence $\{u^{\var,\delta}\}$ 
converges to
an entropy solution $u$ of the first order hyperbolic conservation law
$$
\pa_tu+\pa_xf(u)=0
\tag 1.3
$$
as $\var$ and $\delta\to 0$. 
Throughtout the paper, we assume the following conditions on $f$ and $\beta$: 
\item{(A$_1$)} There are two constants $C_1>0$ and $m>1$ such that 
$$
|f'(u)|\le C_1 \, (1+|u|^{m-1})  \qquad \text { for all } \, u\in{\R}. 
$$
\item{(A$_2$)} The function $\beta=\beta(\lambda)$ is non-decreasing and 
satisfies 
$\beta(0)=0$.

Observe that, from (A$_2$), we can deduce $\beta(\lambda)\lambda\ge0$ for all
$\lambda\in{\R}$.
The following assumption will also be in use:  
\item{(B$_1$)} There exist constants $C_2$, $C_3$, $N>0$, and $r\ge 1$ such
that
$$
C_2 \, |\lambda|^{3r}\le \beta(\lambda)\lambda\le C_3 \, |\lambda|^{3r} 
\qquad \text{ for all } \, |\lambda|\ge N.
$$
The function $\beta$, therefore, is at least
quadratic and could vanish on a bounded interval. 
Two other assumptions on the function  $\beta$ will be 
required for certain results below:
\item{(B$_2$)} There are constants $C_4>0$ and $N>0$ such that
$$
\beta(\lambda)\lambda\ge C_4 \, |\lambda|^3
\qquad \text{ for all } \, |\lambda|\ge N.
$$
\item{(B$_3$)} There are constants $C_5>0$ and $r\ge 1$ such that
$$
\beta(\lambda)\lambda\ge C_5 \, |\lambda|^{3r} 
\qquad \text{ for all } \, \lambda\in{\R}.
$$

Observe that both (B$_2$) and (B$_3$) are refinements to the lower bound 
estimate in
(B$_1$), and that (B$_2$) is actually a consequence of (B$_3$). 
For our main theorems, we shall assume either (B$_1$), or
(B$_1$) and (B$_3$). The condition (B$_2$) will be relevant only in the course of 
the 
discussion of the proofs. 
One typical example of a function $\beta$ 
satisfying (A$_2$) and (B$_1$) (and also (B$_2$), (B$_3$)) is given by
$$
\beta(\lambda)=|\lambda|^{3r-2}\lambda\quad(r\ge1).
\tag 1.4
$$

Recall that, for the Cauchy problem for (1.3),  
existence and uniqueness of an entropy 
solution were first proved by Kruzkov \cite{\refKruzkov} (Section 2 below).
The {\it vanishing viscosity limit\/}, which corresponds to $\delta=0$ and
$\var\to0$, was studied by several authors in the special case of 
a linear (and therefore non-degenerate) diffusion
term $\beta(\lambda)=\lambda$; this activity started with 
Oleinik's \cite{\refOleinik} and Kruzkov's \cite{\refKruzkov} works. 
For more general viscosity coefficients and under the 
assumptions (A$_2$)--(B$_1$),  the convergence was proved in \cite{\refMarcati}.

The {\it vanishing dispersion limit\/}
($\var=0$ and $\delta\to 0$) has also been widely 
investigated, see Lax-Levermore \cite{\refLaxL}, 
\cite{\refDiPerna, \refVenakidesOne, \refVenakidesTwo}, and the survey paper 
by Lax \cite{\refLaxTwo}. 
The main source of difficulty lies 
in the highly oscillating behavior of the solutions 
$u^{0,\delta}$, which do not converge in a strong topology.

The complete equation (1.1) arises, at least for the linear diffusion
$\beta(\lambda)=\lambda$ (the so called Korteweg-de Vries-Burgers equation), 
as a
model of the nonlinear propagation of dispersive and dissipative waves in many
different physical systems. For related 
studies, let us refer to \cite{\refBonaP, \refBonaS, \refCanosa, \refJeffrey, 
\refJohnson}. 
Travelling waves were studied in \cite{\refBonaS} and 
\cite{\refJacobs} and (partial)
analytical results can be found in \cite{\refKenigThree} and \cite{\refCraig}. 
For a study of numerical approximations with similar features as 
(1.1), see \cite{\refChen}.

The vanishing dispersive and
diffusion limit (both $\var$ and $\delta$ tend to zero) was first studied 
by Schonbek \cite{\refSchonbek}. 
She assumed  $\beta(\lambda)=\lambda$ and used the compensated compactness
method introduced by Tartar \cite{\refTartar}.  
In particular, she proved that the sequence $\{u^{\var,\delta}\}$ 
converges to a weak solution to (1.3) 
(but not necessarily the unique
entropy one) under the assumption that either $\delta=0(\var^2)$ for 
$f(u)= u^2 / 2$, or
$\delta=0(\var^3)$ for arbitrary subquadratic flux-functions $f$ 
(i.e.~take $m=2$ in (A$_1$)). Our inequalities on $\var$ and $\delta$ requires
that viscosity dominates dispersion, which is expected since  
$\{u^{\var,\delta}\}$ 
are known not to converge to a weak solution of (1.3) (a~fortiori to the entropy
one)
when dispersion effects are dominant.  

Here we consider the general nonlinear equation (1.1) and we show that, under 
the 
assumption
(B$_1$) and if $\delta\le C \, \var^{5-m\over 3-m}$ $(m<3)$, the solutions
$u^{\var,\delta}$ converge to the entropy weak solution of (1.3) (Theorem 4.1). 
Under the assumptions (B$_1$) and (B$_3$) and 
when $\delta\le C \, \var^{5-m\over r(5-m)-1}$ $(m<5-{1/ r},\ r\ge1)$, 
the same result holds true (Theorem 4.2). Furthermore, 
our Theorem 4.3 below improves upon 
Schonbek's result in \cite{\refSchonbek} in the case that 
$\beta(\lambda)=\lambda$ and  $f$ is subquadratic.

The main tool necessary to deal with these singular limits
is the concept of entropy 
measure-valued solution to the equation (1.3). 
Measure-valued solutions were introduced by DiPerna \cite{\refDiPerna} (see 
Section 2)
to represent the weak-$\star$ limits of the solutions of (1.1),  
According to DiPerna's theory, the strong convergence of an
approximate sequence follows once one obtains
\item{(i)} uniform bounds in $L^\infty(0,T;L^q({\R}))$ for all $T>0$ with $q>m$ 
(where $m$ is given by (A$_1$);
\item{(ii)} weak consistency of the sequence with the entropy inequalities and 
the
strong consistency with the initial data (Cf.~the conditions (2.5) and (2.6) of 
Section 2).

Such a framework can be used to establish 
as well the convergence of difference schemes as shown in \cite{\refCoquel}. 
Observe however that, because of the oscillatory behavior of $u^{\var,\delta}$, 
no
maximum principle is available for the equation (1.1) 
and the $L^q$ (uniform) estimates  are  
more natural. Therefore we actually use the $L^q$ version of DiPerna's result
derived by Szepessy in \cite{\refSzepessy}. 
\medskip

One basic estimate that does not require any restriction on the parameters 
$\var$ and $\delta$ is
the so-called energy estimate based on the simplest quadratic function $u^2$ 
(Lemma 3.1). 
The derivation of the necessary $L^q$ estimates with larger exponents ($q >2$)
is based on the nonlinearity
of the viscosity coefficient $\beta$ and requires suitable 
restriction on $\var$ and $\delta$. 

The problem studied in this paper is motivated by the numerical method 
proposed by Von Neumann and Richtmyer
\cite{\refVonNeumann} (see also \cite{\refRichtmyer}) 
and based on the pseudo-viscosity  approximation
$\beta(\lambda)=|\lambda|\lambda$.
An important advantage of the
nonlinear (quadratic) diffusion term $\beta(\lambda)$ is that it adds sharply 
localized
viscosity  
near shocks  and a small quantity (possibly zero if $\beta$ vanishes for 
$\lambda<0$) 
elsewhere. The
pseudo-viscosity idea  was proposed by Von Neumann 
to stabilize an unstable scheme such as the 
Lax-Wendroff scheme \cite{\refLaxW, \refLaxOne, \refGoodman}: 
a minimal amount of numerical viscosity is added in order to prevent both
the formation of unphysical (entropy violating) shocks and the generation of 
highly
oscillatory approximate solutions. 

\beginsection0{2.}{Preliminaries}

This section contains short background material on Young measures, entropy 
measure 
valued (mv) solutions (Section 2.1), and dispersive equations (Section 2.2).
\bigskip

\noindent{\bf  2.1 Young Measures and Measure-Valued Solutions}

Following Schonbek \cite{\refSchonbek}, 
we describe a representation theorem for Young
measures associated with a sequence of uniformly bounded functions of
$L^q$ where $q \in (1,\infty)$ is fixed in the whole of this subsection.
The corresponding setting in $L^\infty$ was first established by Tartar
\cite{\refTartar}.

\proclaim{Lemma 2.1} Let $\{u_j\}$ be a uniformly bounded sequence in
$L^\infty({\R}_+$; $L^q ({\R}))$. Then there exists a subsequence $\{u_{j'}\}$ and a
weak-$\star$ measurable,
mapping $\nu:{\R}\times{\R}_+\to \t{\rm Prob }({\R})$ taking 
its values in the spaces of non-negative measures
with total mass one (probability measures) such
that, for all functions 
$g\in \Cal C({\R})$ satisfying 
$$
g(u)=0(|u|^r)  \quad \t{\rm as }|u|\to \infty\ 
\tag 2.1
$$
for some $r\in [0,q)$, the following limit representation holds
$$
\lim_{j'\to\infty}\iint_{{\R}\times{\R}_+} g(u_{j'}(x,t))\phi(x,t)\d x\d t
=\iint_{{\R}\times{\R}_+} \int_{\R} g(\lambda)\d \nu_{(x,t)}(\lambda) \, \phi(x,t)\d 
x\d t
\tag 2.2
$$
for all $\phi\in C_0^\infty({\R}\times{\R}_+)$.
\endproclaim

The measure-valued function $\nu_{(x,t)}$ is a Young measure associated with the
sequence $\{u_j\}$. The following result reveals the connection between the 
structure
of $\nu$  and the strong convergence.

\proclaim{Lemma 2.2} Suppose that $\nu$ is a Young measure associated with a
sequence $\{u_j\}$ that is uniformly bounded in $L^\infty({\R}_+;\ L^q ({\R}))$. 
For $u \in L^\infty({\R}_+;\ L^q ({\R}))$, the following
statements are equivalent:
\item{\rm (i)} $\lim_{j\to \infty}u_j=u\quad$ in $L^\infty({\R}_+; 
L^r_{\loc}({\R}))$ for some $r\in
[1,q)$;
\item{\rm (ii)} $\nu_{(x,t)}=\delta_{u(x,t)}$ for almost every 
$(x,t)\in{\R}\times{\R}_+$.
\endproclaim

In (ii) above, the notation $\delta_{u(x,t)}$ is used for the Dirac mass defined
by
$$
\iint_{{\R}\times{\R}_+} <\delta_{u(x,t)}, g(\cdot)>  \, \phi(x,t)\d x\d t
= 
\iint_{{\R}\times{\R}_+}  g(u(x,t))  \, \phi(x,t)\d x\d t
$$
for all $g\in \Cal C({\R})$ satisfying (2.1) and all $\phi\in 
C_0^\infty({\R}\times{\R}_+)$. 
Following DiPerna \cite{\refDiPerna} and \cite{\refSzepessy}, 
we now define the measure-valued solutions to the first order Cauchy problem
$$
\pa_tu+\pa_x \,f(u)=0,
\tag 2.3
$$
$$
u(x,0)=u_0(x), \qquad x \in {\R}.
\tag 2.4
$$

\proclaim{Definition 2.1} 
Assume that $f$ satisfies the growth condition {\rm (2.1)} and 
$u_0\in L^1({\R})\cap L^q ({\R})$. 
A Young measure $\nu$ associated with a sequence $\{u_j\}$, which is assumed to 
be 
uniformly bounded in $L^\infty({\R}_+;L^q ({\R}))$, is called an entropy mv solution
to
the problem {\rm (2.3)-(2.4)} if 
$$
\pa_t\la\nu_{(\cdot)},|\lambda-k|\ra+\pa_x
\la\nu_{(\cdot)},\sgn(\lambda-k)(f(\lambda)-f(k))\ra\le 0
\tag 2.5
$$
in the distributional sense for all $k\in{\R}$ and, for all interval $I\subseteq 
{\R}$,
$$
\lim_{T\to 0^+}{1\over T}\int^T_0\int_I \la \nu_{(x,t)}, |\lambda-u_0(x)|\ra\d 
x\d t=0.
\tag 2.6
$$
\endproclaim

\proclaim{Remark} \rm 
A function $u\in L^\infty({\R}_+;L^1({\R})\cap L^q ({\R}))$ is 
an entropy weak solution to (2.3)-(2.4) in the sense of Kru\v{z}kov 
\cite{\refKruzkov} if and only if the Dirac measure
$\delta_{u(\cdot)}$ is an entropy mv solution. In the case $p=+\infty$, 
existence and
uniqueness of such solutions were proved in \cite{\refKruzkov}. The following 
results of entropy mv solutions were proved in \cite{\refSzepessy}:
Theorem 2.3 states that entropy mv solutions are actually Kruzkov's solutions. 
Theorem 2.4 states that the problem has a unique $L^q$ solution. 
\endproclaim

\proclaim{Theorem 2.3} 
Assume that $f$ satisfies {\rm (2.1)} and $u_0\in L^1({\R})\cap L^q ({\R})$. 
Suppose that $\nu$ is an entropy mv
solution to {\rm (2.3)-(2.4)}.
Then there exists a function $w\in L^\infty({\R}_+;
L^1({\R})\cap L^q ({\R}))$ such that
$$
\nu_{(x,t)}=\delta_{w(x,t)}
\tag 2.7
$$
for almost every $(x,t)\in{\R}\times{\R}_+$.\endproclaim

\proclaim{Theorem 2.4} 
Assume that $f$ satisfies {\rm (2.1)} and $u_0\in L^1({\R})\cap L^q ({\R})$. 
Then there exists a unique entropy solution
$u \in L^\infty({\R}_+;L^1({\R})\cap L^q ({\R}))$
to {\rm (2.3)-(2.4)} which, moreover, satisfies 
$$
\|u(\cdot,t)\|_{L^r({\R})}\le \|u_0\|_{L^r({\R})}
\tag 2.8
$$
for almost every $t\in{\R}_+$ and for all $r\in[1,q]$.
Moreover the measure-valued mapping $\nu_{(x,t)}=\delta_{u(x,t)}$ is 
the unique 
entropy mv solution of the same problem.
\endproclaim

Combining Theorems 2.3 and 2.4 and Lemma 2.2 we obtain our main convergence 
tool.

\proclaim{Corollary 2.5} 
Assume that $f$ satisfies {\rm (2.1)} and $u_0\in L^1({\R})\cap L^q ({\R})$. 
Let $\{u_j\}$ be a sequence of functions that are 
uniformly bounded in
$L^\infty({\R}_+;L^q ({\R}))$ for   $q\geq 1$, 
and let $\nu$ be a Young measure associated with this sequence. 
If $\nu$ is an
entropy mv solution to {\rm (2.3)-(2.4)}, then
$$
\lim_{j\to\infty}u_j=u\quad\t{\rm in }L^\infty({\R}_+;L^r_{loc}({\R}))
\tag 2.9
$$
for all $r\in [1,q)$, where $u \in L^\infty({\R}_+;L^q ({\R}))$
is the unique entropy solution to {\rm (2.3)-(2.4)}.
\endproclaim
\bigskip

\noindent{\bf 2.2 Nonlinear Dispersive Equations.}

Consider the fully nonlinear, KdV-type equation in one space variable 
$$
\pa_tu+F(u,\pa_xu,\pa^2_x u,\pa_x^3 u)=0
\tag 2.10
$$
for $(x,t)\in{\R}\times(0,T)$ together with the Cauchy data
$$
u(x,0) = u_0(x) \qquad \text{ for } \, x \in {\R}. 
\tag 2.11
$$

In recent years this class  of problems has been extensively studied;
see for instance
\cite{\refKenigOne, \refKenigTwo, \refKenigThree}
and the references cited therein. In particular, under suitable
assumptions, these equations enjoy a gain of regularity of the
solutions with respect to their initial data 
(Cf.~\cite{\refCraig, \refKato, \refKruzkovF} and the references above).
We observe however that a complete theory of global existence remains 
to be developped. Presenting a complete review is far 
beyond the aim of this presentation. Here 
we need only recall a result --that is sufficient to deal with (1.1)-- 
of local existence and uniqueness ot the smooth solutions to
the equation (2.10). 
More details and proofs can be found in the paper by 
Craig-Kappeler-Strauss \cite{\refCraig} 
(see also \cite{\refKenigThree}).

Using the notation $U=(u_0,u_1,u_2,u_3)$, 
the assumptions on the smooth function $F=F(U)$
are as follows:
\item{(H$_1$)} There exists $\delta>0$ such that
$$
\pa_{u_3}F(U)\ge \delta>0 \qquad \text{ for all } \, U\in {\R}^4; 
$$
\item{(H$_2$)} $\pa_{u_2}F(U)\le 0 \qquad $ for all $U\in{\R}^4$.

\noindent The equation (1.1) does satisfy these hypotheses provided
$\delta>0$ and (A$_2$) holds.

\proclaim{Theorem 2.6 {\rm (Uniqueness)}} Let $T>0$ be fixed and
assume $F$ satisfies {\rm (H$_1$)-(H$_2$)}. 
For any $u_0\in H^7({\R})$ there is at most one solution
$u\in L^\infty((0,T);H^7({\R}))$ to the problem {\rm (2.10)-(2.11)}.
\endproclaim

\proclaim{Theorem 2.7 {\rm (Existence)}} 
Assume $F$ satisfies {\rm (H$_1$)-(H$_2$)}.
Let
$N$ be an integer $\ge7$ and let $C_0>0$ be a given constant. 
There exists a time $T>0$,  depending
only on $C_0$, such that for all $u_0\in H^N({\R})$ with $\|u_0\|_{H^7}\le C_0$, 
there
exists at least one solution $u\in L^\infty((0,T);H^N({\R}))$ to the problem 
{\rm (2.10)-(2.11)}.
\endproclaim

Observe that the setting above based on a Sobolev norm 
of relatively high order ($H^7({\R})$) is the optimal result provided by 
the current techniques of analysis of dispersive equations. 
The {\it global\/} existence of smooth solutions to (1.1)-(1.2) appears to be an
open problem. In the following we tacitly 
restrict attention to a time $T_* \in (0, \infty]$ chosen such that the problem 
(1.1)-(1.2)
is well-posed in the strip ${\R}\times(0,T_*)$.

\beginsection0{3.}{A~Priori Estimates}

In this section we consider a sequence $\{u^{\var,\delta}\}$ of smooth solutions
of
(1.1)-(1.2) that vanish at infinity. We assume also that the initial data
$\{u_0^{\var,\delta}\}$ are smooth functions with compact support, and are 
uniformly
bounded in $L^1({\R})\cap L^q ({\R})$ for a suitable $q>1$.
In what follows, whenever it does not lead to
confusion, we omit the indices $\var$ and $\delta$. Similarly all
constants are denoted by $C, C_1, \cdots$. 
We begin with the natural energy estimate 
based the quadratic function $u^2/2$.  (Arbitrary convex functions 
could not be easily used here because of the dispersive term.)

\proclaim{Lemma 3.1} For any $T>0$, it holds that 
$$
\int_{\R} u^2(x,T)\d x + 2 \, \var\int^T_0\int_{\R}\beta(u_x(x,t))u_x(x,t)\d x\d 
t  = \int_{\R} u^2_0(x)\d x.
\tag 3.1
$$
\endproclaim

\noindent{\it Proof.\/} 
We multiply (1.1) by $u$ and integrate in space. Integrating by
parts we obtain
$$
{\d\over\d t}\ {1\over2}\int_{\R} u^2\d x+\int_{\R} uf'(u)u_x\d x
=-\var\int_{\R}\beta(u_x)u_x\d x-\delta\int_{\R} u\,u_{xxx}\d x.
\tag 3.2 
$$
The second terms on both sides of (3.2) vanish identically, 
since we can write these terms in a conservative form:
$$u\,f'(u)u_x=(G(u))_x
$$
with $G'=u\,f'(u)$ and 
$$
u\,u_{xxx}=(u\,u_{xx}-{1\over2}u^2_x)_x.
$$
The estimate (3.1) follows therefore from (3.2).
{$\square$}\medskip

At this stage, in view of the assumption (A$_1$) and Lemma 3.1,
we have a uniform control of $u$ in
$L^\infty({\R}_+;L^2({\R}))$ and $\var \, \beta(u_x)u_x$ in $L^1((0,T)\times {\R})$
for every $T>0$.
Our next aim is to derive a (non-uniform) 
estimate of $u$ in $L^\infty$ norm, which will be uniform in
$\var$ but not in $\delta$. We first provide a second energy-type estimate.

\proclaim{Lemma 3.2} Let $F$ be defined by $F'(u)=f(u)$. 
For every $T>0$, we have
$$
\aligned
&{\delta\over 2} \int_{\R} u^2_x(x,T)\d x-\int_{\R} F(u(x,T))\d x +\var \, 
\delta\int^T_0\int_{\R}
u^2_{xx}\beta'(u_x)\d x\d t\\
&\quad ={\delta\over 2} \int_{\R} u^2_{0,x}(x)\d x-\int_{\R} F(u_0(x))\d x+\var 
\int^T_0\int_{\R}
f'(u)\beta(u_x)u_x \, \d x\d t. 
\endaligned
\tag 3.3
$$
\endproclaim

\noindent{\it Proof.\/} Multiplying (1.1) by $f(u)+\delta u_{xx}$ we obtain the
equality
$$
\align
0&=\int_{\R}\{F(u)_t+H(u)_x-\var f(u)\beta(u_x)_x+\delta f(u) u_{xxx}\}\d x\\
&+\int_{\R}\{\delta\, u_{xx}u_t+\delta u_{xx} f(u)_x-\delta 
u^2_{xx}\,\beta'(u_x)+\delta
u_{xxx}u_{xx}\}\d x
\endalign
$$
where we have set $H'=ff'$. After integration by parts in space, we have
$$
{\d\over\d t}\int F(u) \d x 
-{\delta\over 2} \, u^2_x\,\d x+\int_{\R}\var  f'(u)\beta(u_x)u_x\, \d x
- \var \, \delta \int_{\R} u^2_{xx}\beta'(u_x)\d x=0,
$$
which gives the desired conclusion.
{$\square$}\medskip

Using the bound for $\sqrt{\delta}\;u_x$ in $L^\infty({\R}_+;L^2({\R}))$ 
that follows from Lemma 3.2, 
we are now able to estimate $u$ in the $L^\infty$ norm.

\proclaim{Lemma 3.3} If $m<5$ in  the assumption {\rm (A$_1$)}, there exists a
constant $C>0$ such that 
$$
|u(x,t)|\le C \, \delta ^{-{1\over 5-m}}
\tag3.4
$$
for all $(x, t)\in {\R}\times{\R}_+$.
\endproclaim

\noindent{\it Proof.\/} In view of (A$_2$) and for all $u\in {\R}$, 
we have
$$
|F(u)|\le C \, (1+|u|^{m+1}).
$$
The main idea is to use (3.3) to control $\delta \, u^2_x$ in terms of 
$F(u)$, the latter being estimated from the above growth condition. 
We deduce from (3.3) that 
$$
\align
&{\delta\over 2} \int_{\R} u^2_x(x,T)\d x+\var \, \delta\int^T_0\int_{\R} 
u^2_{xx}\,\beta'(u_x)\d
x \d t\\
&\qquad ={\delta\over 2} \int_{\R} u^2_{0,x}(x)\d x+\int_{\R} 
\big(F(u(x,T))-F(u_0(x)) \big)\d x
+\var\int^T_0\int f'(u)\beta(u_x)u_x\,\d x\d t\\
&\qquad\le C \, + \, C \, \|u(\cdot ,T)\|^{m-1}_{L^\infty({\R})}\left(\|u(\cdot
,T)\|^2_{L^2({\R})}+\var\int^T_0\int_{\R}\beta(u_x)u_x\,\d x\d t\right), 
\endalign
$$
where we used (A$_2$) as well. 
(We do not explicitly write the terms involving
the initial data
as it is assumed to be a smooth function.)
Therefore in view of 
(3.1), we arrive at
$$
{\delta\over 2} \int_{\R} u^2_x(x,T)\d x 
+ \var \, \delta\int^T_0\int_{\R} u^2_{xx}\,\beta'(u_x)\d x\d t
\le C \, + C \, \|u(\cdot, T)\|^{m-1}_{L^\infty({\R})}.
$$
Using the Cauchy-Schwartz inequality, we obtain
$$
\align
|u(x,t)|^2&\le 2\int^x_{-\infty}|u(y,t)u_x(y,t)|\d y\\
&\le {2\over\sqrt{\delta}}\|u(\cdot, t)\|_{L^2({\R})}\sqrt{\delta}\,\|u_x(\cdot, 
t)\|_{L^2({\R})}\\
&\le {C\over \sqrt{\delta}}\|u_0\|_{L^2({\R})}
\left(1+\|u(\cdot,t)\|^{m-1}_{L^2({\R})}\right)^{1 / 2}.
\endalign
$$
Hence, for all $t>0$,
$$
\|u(\cdot,t)\|^4_{L^\infty({\R})}\le{C\over \delta}\left(1+\|u(\cdot,
t)\|^{m-1}_{L^\infty({\R})}\right).
\tag 3.5
$$

Now, since $m<5$, the growth of the left hand side of (3.5) 
exceeds the growth of the right hand side. So we
 can check that (3.5) implies 
$$
\align
\|u(\cdot,t)\|_{L^\infty({\R})}&\le 
\max\left(1,\left({2C\over\delta}\right)^{1\over
5-m}\right)\\
&\le C\,\delta^{-{1\over 5-m}}.\endalign
$$
Namely, setting $y=\|u(\cdot,t)\|_{L^\infty}$, we have  for $y>0$
$$
y^4\le{C\over\delta}\left(1+y^{m-1}\right).
\tag 3.6
$$
If $y\le 1$ we have the conclusion. Otherwise, suppose we had 
$y>\left({2C\over\delta}\right)^{1\over 5-m}$, then we would deduce
that
$$
y^4>{2C\over \delta}\,y^{m-1}>{C\over\delta}\left(1+y^{m-1}\right),
$$
which would contradict the inequality (3.6).
{$\square$}\medskip

In view of the proof of Lemma 3.3 we also state the following result:

\proclaim{Lemma 3.4} For any $T>0$ we have
$$
{1\over 2}\int_{\R} u^2_x(x,T)\d x+\var \int^T_0\int_{\R} u^2_{xx}\beta'(u_x)\d x\,
\d t \le
C \, \delta^{-{4\over 5-m}}.
\tag 3.7
$$
\endproclaim

Finally we derive uniform bounds in $L^\infty((0,T)$; $L^q({\R}))$ 
with $q\le 5$. The following result provides a uniform estimate in the 
Lebesgue space $L^5({\R})$ (so improving upon the $L^2$ bound in Lemma 3.1) 
by taking advantage of the nonlinearity property of $\beta$ 
(for large values of $\lambda=u_x$) as stated in {\rm (B$_2$)}.

\proclaim{Proposition 3.5} Assume that the assumption {\rm (B$_2$)} holds
and $m<3$ in the assumption {\rm (A$_1$)}. 
There exists a constant $C>0$, which depends only on the initial data, 
such that, for all $\delta$ and $\var$ small enough and for every $T>0$, we have
$$
\sup_{t\in(0,T)}\|u(\cdot, t)\|^5_{L^5({\R})}\le C \, \left( T 
+\var^{-1}\delta^{3-m\over
5-m}\right).
\tag 3.8
$$
\endproclaim

\noindent{\it Proof.\/} 
Set $\eta(u)=|u|^5$. Multiply (1.1) by $\eta'(u)$ and integrate
on ${\R}\times(0,T)$. It follows easily that 
$$
\aligned
&\int_{\R}\eta(u(x,T))\d x+\var \int^T_0\int_{\R}\eta''(u)\beta(u_x)u_x\,\d x\d t \\
&\qquad=\int_{\R}\eta(u_0(x))\d x-{\delta\over 
2}\int^T_0\int_{\R}\eta'''(u)u_x^3\,\d x\d t.
\endaligned
\tag 3.9
$$
On the other hand, according to (B$_2$), 
there exists a constant $C$ such that for all $\lambda \in {\R}$ 
$$
|\lambda|^3\le C \, (1+\beta(\lambda)\lambda).
$$
Using the latter in the energy estimate (3.1), we are able to control
the second term in the right hand side of (3.9):
$$
\align
\left|\int^T_0\int\eta'''(u)u^3_x\,\d x\d t\right|
& \le C_1 \, \int^T_0\int_{\R} |u|^2 |u_x|^3 \d x\d t \\
& \le C_2 \, \int^T_0\int_{\R} |u|^2 \big( 1 + \beta(u_x)u_x \big) \,\d x\d t \\
&\le C_3 \, T \sup_{t \in (0,T)} \|u(t)\|_{L^2({\R})}
           + {C \over \var}\|u\|^2_{L^\infty({\R}\times(0,T))}\\
& \le C \, T \, + C \, \left(1+\var^{-1}\delta^{3-m\over 5-m}\right).
\endalign
$$
where the latter inequality follows from the $L^\infty$ bound in Lemma 3.3.
Returning to (3.9), we obtain (3.8).
{$\square$}\medskip

Under the stronger assumption (B$_3$), we have a sharper estimate.

\proclaim{Proposition 3.6} Assume that {\rm (B$_3$)} holds 
for a given $r\ge1$ and 
that $m<q \equiv 5-{1/ r}$ in {\rm (A$_1$)}. There exists a 
constant $C>0$, 
which depends only on the initial data, such that, for all $\delta$ and $\var$ 
small enough and every $T>0$, 
$$
\sup_{t\in(0,T)}\|u(\cdot, t)\|_{L^q({\R})}^q
     \le C \, \left(T^{1-{1 / r}} 
         + \var^{-{1/ r}}\delta^{-{q-m \over 5-m}}\right).
\tag 3.10
$$
\endproclaim

\noindent{\it Proof.\/} As in the proof of Proposition 3.5, 
we consider the formula (3.9)
but now with $\eta(u)=|u|^q$ with $q=5-{1/ r}$. 
The assumption (B$_3$) yields, using (3.1) and (3.4),
$$
\align
\left|\int^T_0\int_{\R}\eta'''(u)u^3_x\,\d x\d t\right| 
& \le 
      C\int^T_0\int_{\R}|u|^{q-3}|u_x|^3\d x\,\d t\\
& \le C \, \left(\int^T_0\int_{\R} |u|^{r'(q-3)}\d x\d t\right)^{1 / r'}
           \left(\int^T_0\int_{\R} \beta(u_x)u_x\,\d x\d t\right)^{1 / r}\\
& \le C \, \var^{-{1/ r}}\|u\|^{1/ r}_{L^\infty({\R})} \, T^{1/r'} \, 
      \sup_{t\in(0,T)}\|u(\cdot,t)\|_{L^2({\R})}^{2 /  r}\\
& \le C \, T^{1/r'} \, \var^{-{1 / r}}\delta^{-{1\over r(5-m)}}, 
\endalign
$$
where $1/r + 1/r' = 1$ and since $r'(q-3) = 2 + r'/r$. Now (3.10) follows 
from the formula (3.9). 
{$\square$}\medskip

Observe that the estimates (3.8) and (3.10) hold only thanks to the nonlinear 
form of the viscosity term. Note that when $\delta=0$, these
estimates reduce to bound that are uniform in $\var$, which 
is expected for conservation laws with viscosity but no dispersion. On the other
hand, 
most of our estimates blow-up when taking $\var=0$ and 
provide no control of $L^q$ norms. Observe however that the $L^2$ bound in Lemma
3.1
is uniform for $\delta$ arbitrary and $\var \to 0$. 
Loosely speaking, this is also the case
of the estimate in Proposition 3.6 if $r$ could be
chosen to be $r = \infty$, in which case (3.10) becomes a uniform $L^5$ 
estimate. 
Such a choice of $r$ is not allowed however, cf.~(B$_3$).

For the sake of completeness, we finally state an 
analogous estimate for linear diffusions
which was proved by Schonbek \cite{\refSchonbek}.

\proclaim{Proposition 3.7} Let $\beta(\lambda)=\lambda$ and take $m=2$ in {\rm
(A$_1$)}. For any $T>0$, there exists a constant $C_T>0$, which depends 
only on the initial data, such that for $\delta\le \nu\,\var^3$
$$
{\sup_{t\in(0,T)}} \|u(\cdot,t)\|_{L^4({\R})}\le C_T.
\tag 3.11
$$
\endproclaim
\bigskip
\beginsection0{4.}{Convergence Results}

In Section 3, we have established several uniform bounds for the sequence
$\{u^{\var,\delta}\}$ of solutions to the Cauchy problem (1.1)-(1.2) under 
certain assumptions on the functions $f$, $\beta$, and the parameters
$\var$ and $\delta$.
Assume again that the initial data $u_0^{\var,\delta}$
are smooth with compact support and
that there exists a limiting function $u_0\in L^1({\R})\cap L^q ({\R})$ and a 
suitable 
$q>1$ (specified below) such that, if $\delta=O(\var)$,  
$$\lim_{\var\to 0} u_0^{\var,\delta}= u_0 
\quad  \text{ in }L^1({\R})\cap L^q ({\R}).
\tag 4.1
$$ 
Returning to the proofs of Section 3, it is not hard to see that the following
conditions on the initial data are sufficient for the estimates therein to hold
uniformly 
with respect to a class of initial data:
$$
\|u_0^{\var,\delta}\|_{L^2({\R})} \, + \, \|u_0^{\var,\delta}\|_{L^q({\R})} 
\delta^{1/2} \, \|u_{0,x}^{\var,\delta}\|_{L^2({\R})}  \, \le\, C. 
$$
In this section we prove the strong convergence 
of the sequence $u^{\var,\delta}$.

\proclaim{Theorem 4.1} Assume that {\rm (B$_1$)} holds and $m<3$ in {\rm 
(A$_1$)}. Let
$u^{\var,\delta}$ be a sequence of smooth solutions to {\rm (1.1)-(1.2)} on
${\R}\times(0,T)$ (for a given $T>0$), 
which vanish at infinity and are associated with initial data satisfying
{\rm (4.1)} with $q=5$.
If there is a constant $C>0$ such that $\delta\le C \, \var^{5-m\over 3-m}$, 
then the (whole) sequence $u^{\var,\delta}$ converges to a function 
$u \in L^\infty((0,T);L^5({\R}))$,
which is the unique entropy solution to {\rm (2.3)-(2.4)}.
\endproclaim

\proclaim{Theorem 4.2} Assume that {\rm (B$_1$)} and {\rm (B$_3$)} hold 
and $m<5-{1/  r}$ in {\rm (A$_1$)}.
Let $u^{\var,\delta}$ be a sequence of smooth solutions to {\rm (1.1)-(1.2)} on
${\R}\times(0,T)$ (for a given $T>0$),
which vanish at infinity and are associated with 
initial data satisfying
{\rm (4.1)} with $q=5-{1/  r}$. 
If there is a constant $C>0$ such that 
$\delta\le C \, \var^{5-m\over r(5-m)-1}$, 
then the (whole) sequence converges to a function
$u\in L^\infty((0,T);L^q({\R}))$, $q=5-{1/ r}$, which is the unique entropy 
solution
to {\rm (2.3)-(2.4)}. \endproclaim

Let us give an analogous statement in the case $\beta(\lambda)=\lambda$, which
improve upon \cite{\refSchonbek} (Cf.~Theorem 5.1 therein).

\proclaim{Theorem 4.3} Assume $m=2$ in {\rm (A$_1$)} and let
$\beta(\lambda)=\lambda$. If 
$\delta\le C \, \var^3$, the whole sequence
$u^{\var,\delta}$ of solutions of {\rm (1.1)-(1.2)} converges to a function 
$u\in
L^\infty((0,T);L^4({\R}))$, 
which is the unique entropy solution to {\rm (2.3)-(2.4)}.
\endproclaim

Similar convergence results can be proven for the case $m=3$ or for some special
flux-functions along the lines of what was done 
in \cite{\refSchonbek} (Cf.~Sections 4 and 5 therein).
Theorems 4.1-4.2 follow easily by simply using the following general result 
and the $L^q$ bounds derived in Section 3, 
Proposition 3.5 and Proposition 3.6 respectively.

\proclaim{Theorem 4.4} Assume that {\rm (B$_1$)} holds. Let
$u^{\var,\delta}$ be a sequence of smooth solutions to {\rm (1.1)-(1.2)} on
${\R}\times(0,T)$ (for a given $T>0$)
associated with initial  data satisfying {\rm (4.1)}. If the sequence is
uniformly bounded in
$L^\infty((0,T);L^q({\R}))$ for $q>m$ and $\delta=o(\var^{1/ r})$, then 
the (whole) sequence converges to a function  $u\in L^\infty((0,T);L^q({\R}))$,
which is the unique entropy solution to {\rm (2.3)-(2.4)}.
\endproclaim

\noindent{\it Proof of Theorem 4.4.\/} 
First of all let us establish that, for any convex function $\eta=\eta(u)$ 
such that $\eta'$, $\eta''$, $\eta'''$ are uniformly bounded on ${\R}$, we have 
$$
\Lambda^{\var,\delta}=\pa_t\eta(u^{\var,\delta})+\pa_x
Q(u^{\var,\delta})\rightharpoonup 0 \quad\t{\rm in } {\Cal D'}({\R}\times{\R}_+),
\tag 4.2
$$
where $Q'=f'\eta'$.
To begin with, observe that
$$
\align
\Lambda^{\var,\delta}&=\var(\eta'(u)\beta(u_x))_x-\var
\eta''(u)\beta(u_x)u_x
-\delta(\eta'(u)u_{xx})_x+\delta\eta''(u)u_{x}u_{xx}\\
&=T_1+T_2+T_3+T_4.
\endalign
$$

To estimate $T_1$, we use the assumption (B$_1$), which implies 
$|\beta(\lambda)|\le C \, \big( 1 + |\lambda|^{3r-1}\big)$ for all $\lambda$. 
For any given $\theta\in C_0^\infty({\R}\times(0,T))$, $\theta\ge 0$, and for
$p={3r\over 3r-1}>1$ ($p'$ being the conjugate exponent of $p$), we get  
$$
\align
\la T_1,\theta\ra &=\left|\iint\var\,\theta_x\eta'(u)\beta(u_x)\d x\,\d
t\right|\\ &\le
C\,\var\|\theta_x\|_{L^{1}({\R} \times (0,T))} + 
C\,\var\|\theta_x\|_{L^{p'}({\R} \times 
(0,T))}\left(\iint|u_x|^{p(3r-1)}\right)^{1/ p},
\endalign
$$
so using the derivative estimate in (3.1) and (B$_1$) again: 
$$
\la T_1,\theta\ra 
\le C \, \big(\var + \var^{1/ 3r}\big)  \, \le 
C \, \var^{1/ 3r}.
\tag 4.3
$$
The second term $T_2$ is nonpositive, namely   
$$
\la T_2,\theta\ra=-\iint\var\eta''(u)\beta(u_x)u_x\theta\ \d x\d t\le 0.
\tag 4.4
$$
To estimate $T_3$, we use the energy estimate (3.1) and the fact that 
$\beta$ is at least quadratic. We write 
$$
\align
\la T_3,\theta\ra&=\delta\iint_{{\R} \times (0,T)} \theta_x\eta'(u)u_{xx}\,\d x\d 
t\\
& =    \delta\iint_{{\R} \times (0,T)}\theta_x\left((\eta'(u)u_x)_x - 
\eta''(u)u^2_x\right) \d x\d t\\
&\le   -\delta\iint_{{\R} \times (0,T)}\theta_{xx}\eta'(u)u_x     
       +    \delta\iint_{{\R} \times (0,T)} |\theta_x| \, u^2_x  \d x\d t  \\
&\le   C\,\delta + \, C \, \delta \, \left(\iint_{\supp\theta}|u_x|^{3r}\d x\d 
t\right)^{1/3r}
         \, + \, C \, \var^{-2/3r}, 
\endalign
$$
where $\supp \theta$ denotes the support of the function in ${\R} \times (0,T)$, 
thus 
$$
\la T_3,\theta\ra 
\le C\,\delta \left(1+C \, \var^{-2/3r} + \var^{-{1/ 3r}}\right) 
\le C \, \left(1+C \, \var^{-2/3r}\right).
\tag 4.5
$$
Finally we deal with $T_4$ as follows: 
$$
\align
\left|\la T_4,\theta\ra\right|&=\left|\delta\iint\eta''(u)u_x\,u_{xx}\theta\,\d 
x\,\d
t\right|\\ 
&=\left|\delta\iint\theta_x
\eta''(u)u^2_x+\delta\iint\theta\eta'''(u){1 \over 2}u^3_x\right|\\
&\le C\,\delta\left(\iint_{\supp\theta}|u_x|^2+|u_x|^3\d x\d t\right),
\endalign
$$
so 
$$
\left|\la T_4,\theta\ra\right| 
\le C\,\delta\left(\var^{-{2/ 3r}}+\var^{-{1/ r}}\right)
\le C \, \var^{-{1/ r}}.
\tag 4.6
$$

Therefore, if $ \delta=o(\var^{1/ r})$, (4.2) follows immediately from
the estimate (4.3)-(4.6). To apply Corollary 2.5 we have to show that (2.5) and 
(2.6) are
satisfied for a Young measure $\nu$ associated with the sequence 
$u^{\var,\delta}$. 
It is a standard matter to deduce, for all convex entropy pairs, 
$$
\pa_t\la\nu_{(\cdot)}, \eta(\lambda) \ra+\pa_x
\la\nu_{(\cdot)}, Q(\lambda) \ra\le 0
$$
from the convergence property (4.2). 
The inequality (2.5) for all $k\in{\R}$ then follows by using a standard
regularization of the function $|u-k|$. Concerning the initial data and in order

to establish (2.6), 
we now combine the entropy inequalities
and the weak consistency property as was suggested by DiPerna 
\cite{\refDiPerna}. We follow the detailled arguments given in 
\cite{\refSzepessy}.

Consider the function 
$g(\lambda)=|\lambda|^r$ for $1<r<\min\t{\rm (2,q)}$, and set
$$
\aligned
G(\lambda,\lambda_0)
& \equiv g(\lambda)-g(\lambda_0)-g'(\lambda_0)(\lambda-\lambda_0)\\
&\ge {r(r-1)\over 2}\ {(\lambda-\lambda_0)^2\over 
(1+|\lambda|+|\lambda_0|)^{2-r}}.
\endaligned
\tag 4.7 
$$
Let $I\subseteq {\R}$ be be a closed and bounded interval. 
Using the Jensen inequality, the Cauchy-Schwartz inequality, 
and the above convexity inequality (4.7), it is easily checked that 
$$
\aligned
&{1\over T}\int^T_0\int_I\la\nu_{(x,t)},|\lambda-u_0(x)|\ra\d x\d t\\
&\qquad \le C_I \, \left({1\over T}\int^T_0\int_I\la 
\nu_{(x,t)},G(\lambda,u_0(x))\ra\d
x\,\d t\right)^{1/ 2}.
\endaligned
\tag 4.8 
$$
Let $\{\psi_n\}\in C_0^\infty({\R})$ be a sequence of test-functions such that
$$
\lim_{u\to\infty}\psi_n = g'(u_0) = \quad \t{\rm in }L^{r'}({\R}),
$$
where $1=1/r + 1/r'$. 
Using the uniform bound in $L^q$ available for the sequence 
$\{u_0^{\var,\delta}\}$, we
get
$$
\aligned
&\int^T_0\int_I \la \nu_{(x,t)},G(\lambda,u_0(x))\ra\d x\d t\\
&\le \int^T_0\int_{\R} \la\nu_{(x,t)}, u_0-\lambda\ra\psi_n\ \d x\,\d
           t+T\int_{{\R}\smallsetminus I}|u_0|^2\d x \\
& + 2 \, T \, \|u_0\|_{L^r({\R})}\|g'(u_0)-\psi_n\|_{L^{r'}({\R})}.
\endaligned
\tag 4.9 
$$
Taking an increasing sequence of compact sets $K_i$ covering ${\R}$, i.e.~such 
that 
$I \subset K_1\subset K_2\subset\dots$ and $\bigcup^\infty_{i=1}K_i={\R}$, we have
$$
\int^T_0\int_I \la \nu_{(x,t)},G(\lambda,u_0(x))\ra\d x\d t
\le\int^T_0\int_{K_i}\la \nu_{(x,t)},G(\lambda,u_0(x))\ra\d x\d t,
$$
which, together with (4.9) where $I$ is replaced by $K_i$, yields
$$
\aligned
& {1\over T}\int^T_0\int_I \la \nu_{(x,t)},Q(\lambda,u_0(x))\ra\d x\d t\\
&\qquad\le{1\over T}\int^T_0\int_{\R} \la \nu_{(x,t)},u_0(x)-\lambda\ra\psi_n\,\d 
x\d t
+2 \, \|u_0\|_{L^r({\R})}\|g'(u_0)-\psi_n\|_{L^{r'}({\R})},
\endaligned
\tag 4.10 
$$
since
$$
\lim_{i \to \infty} \int_{{\R}\smallsetminus K_i}|u_0|^2\d x \, = \, 0. 
$$

Therefore, in view of (4.8) and (4.10), 
the strong consistency property (2.6) will be established if we show that
$$
\lim_{T\to 0^+}{1\over T}\int^T_0\int_{\R}\la 
\nu_{(x,t)},u_0(x)-\lambda\ra\psi_n\,\d x\,\d
t\le 0
\tag 4.11 
$$
for all $n\in \N$. By definition of the Young measure (Cf.~(2.2)), we have  
$$
{1\over T}\int^T_0\int_{\R}\la \nu_{(x,t)},u_0(x)-\lambda\ra \, \psi_n\ \d x\d t
=\lim_{\var, \delta\to 0}{1\over
T}\int^T_0\int \big(u_0(x)-u^{\var,\delta}(x,t)\big) \, \psi_n\, \d x\d t.
$$
On the other hand, we can write 
$$\align
&=\int_{{\R}} \big(u_0(x)-u_0^{\var,\delta}(x,t)\big) \, \psi_n(x)\d x
-{1\over T}\int^T_0\int_{\R}\left(\int^t_0\pa_s\,u^{\var,\delta}(x,s)\d
s\right)\psi_n(x)\d x\d t\\
& \equiv A+B.
\endalign
$$
The term $A$ tends to zero as $\var\to 0$ in view of the weak consistency 
property (4.1). 
Furthermore, by arguing as in the derivation of (4.3), we have 
$$
\align
B&=-{1\over T}\int^T_0\int_{\R}\left(\int^t_0(-\pa_x
f(u^{\var,\delta})+\var\pa_x(\beta(u_x^{\var,\delta}))-\delta\pa^3_x
u^{\var,\delta}\right)\psi_n(x)\d x\d t\\
&=-{1\over
T}\int^T_0\int_{\R} \int^t_0\left(f(u^{\var,\delta})\pa_x\psi_n-
\var\beta(u_x^{\var,\delta})\pa_x\psi_n+\delta
u^{\var,\delta}\pa_x^3\psi_n\right)\d s\,\d x\d t\\
&\le C_n\,T.
\endalign
$$
This leads to the inequality (4.11). The proof of Theorem 4.4 is completed. 
{$\square$}
\medskip

The proof of Theorem 4.3 is based on slightly modified estimates 
in the inequalities
(4.5) and (4.6), which can be derived by arguing as in \cite{\refSchonbek}. 
The details of the proof are omitted.


\bigskip
\noindent{\bf Acknowledgments.} The authors would like to thank Luis Vega for 
helpful conversations about nonlinear dispersive equations, and Brian T. Hayes 
for general discussions on the content of this paper. 
\medskip
\centerline{\smc References}
\medskip\baselineskip=12pt

\ref\key{\refBonaP}
\by \, \,   J.L. Bona, W.G. Pritchard, and L.R. Scott
\paper An evaluation of a model equation for water waves
\jour Philos. Trans. Roy. Soc. London Ser. A
\vol 302
\yr 1981
\pages 457--510
\endref
\ref\key{\refBonaS} 
\by \, \,   J.L. Bona and M.E. Schonbek
\paper Travelling waves solutions to the Korteweg-de Vries-Burgers equation
\jour Proc. Roy. Soc. Edinburgh, Sect. A
\vol 101
\yr 1985
\pages 207--226
\endref
\ref\key{\refCanosa}
\by \, \,   J. Canosa and J. Gazdag
\paper The Korteweg-de Vries-Burgers equations
\jour J. Comput. Phys.
\vol 23
\yr 1977
\pages 393--403
\endref
\ref\key{\refChen} 
\by \, \,   G.Q. Chen and J.G. Liu 
\paper Convergence of difference schemes 
with high resolution to conservation laws
\jour Preprint 
\yr 1994
\endref
\ref\key{\refCoquel} 
\by \, \,   F. Coquel and P.G. LeFloch 
\paper Convergence of finite difference schemes for scalar
conservation laws in several space variables: the corrected
antidiffusive-flux approach
\yr 1991
\vol 57 
\jour  Math. of Comp. 
\pages 169--210 
\endref
\ref\key{\refCraig} 
\by \, \,   W. Craig, T. Kappeler, and W. Strauss
\paper Gain of regularity for equations of KdV type
\jour Ann. Inst. Henri Poincar\'e, Nonlin. Anal. 
\vol 9
\yr 1992
\pages 147--186
\endref
\ref\key{\refDiPerna}
\by \, \,   R.J. DiPerna
\paper Measure-valued solutions to conservation laws
\jour Arch. Rat. Mech. Anal.
\vol 88
\yr 1985
\pages 223--270
\endref
\ref\key{\refGoodman} 
\by \, \,   J. Goodman and P.D. Lax
\paper Dispersive difference schemes I
\jour Comm. Pure Appl. Math.
\vol 41
\yr 1988
\pages 591--613
\endref
\ref\key{\refJacobs} 
\by \, \,   D. Jacobs, B. McKinney, and M. Shearer
\paper Travelling wave  solutions of the modified 
Korteweg-deVries-Burgers equation
\jour J. Diff. Equa.
\vol 116
\yr 1995
\pages 448--467
\endref
\ref\key{\refJeffrey} 
\by \, \,   A. Jeffrey and T. Kakutani
\paper Weak nonlinear dispersive waves: a discussion centered around the 
Korteweg-de Vries equation
\jour SIAM Review
\vol 14
\yr 1972
\pages 582--643
\endref
\ref\key{\refJohnson}
\by \, \,   R.S. Johnson
\paper A nonlinear equation incorporating damping and dispersion
\jour J. Fluid Mech.
\vol 65
\yr 1970
\pages 49--60
\endref
\ref\key{\refKato}
\by \, \,   T. Kato
\paper On the Cauchy problem for the (generalized) 
Korteweg-de Vries equation
\jour Adv. in Math. Suppl. Studies: Studies in Appl. Math.
\vol 8
\yr 1983
\pages 93--128
\endref
\ref\key{\refKenigOne} 
\by \, \,   C.E. Kenig, G. Ponce, and L.  Vega
\paper On the (generalized) Korteweg-de Vries equation
\jour  Duke Math. J.
\vol 59
\yr 1989
\pages 585--610
\endref
\ref\key{\refKenigTwo}
\by \, \,   C.E. Kenig, G. Ponce and L.  Vega
\paper Oscillatory integrals and regularity of dispersive equations
\jour  Indiana  Univ. Math. J.
\vol 40
\yr 1991
\pages 33--69
\endref
\ref\key{\refKenigThree}
\by \, \,   C.E. Kenig, G. Ponce and L.  Vega
\paper Higher-order nonlinear dispersive equations
\jour  Proc. Amer. Math. Soc.
\vol 122
\yr 1994
\pages 157--166
\endref
\ref\key{\refKruzkov}
\by \, \,  S.N. Kru\v zkov
\paper First order quasilinear equations in several independent variables
\jour Mat. Sb.
\vol81
\yr1970
\pages285--255; {\it Math. USSR Sb.} {10} (1970), 217--243
\endref
\ref\key{\refKruzkovF} 
\by \, \,  S.N. Kru\v zkov and A.V. Framinskii
\paper Generalized solutions of the Cauchy problem for the Korteweg-de Vries 
equation
\jour Math. USSR Sb.\vol 48 \yr  1970 \pages 93--131  
\endref
\ref\key{\refLaxOne} 
\by \, \,   P.D. Lax 
\paper On dispersive difference schemes
\jour  Physica D
\yr 1986
\pages  250--254
\endref
\ref\key{\refLaxTwo} 
\by \, \,   P.D. Lax 
\paper The zero dispersion limit, a deterministic analogue of turbulence
\jour  Comm. Pure Appl. Math. 
\vol  44
\yr 1991
\pages  1047--1056
\endref
\ref\key{\refLaxL} 
\by \, \,   P.D. Lax and C.D. Levermore
\paper The small dispersion limit of the Korteweg-de Vries equation
\jour  Comm. Pure Appl. Math. 
\vol 36
\yr 1983
\pages I, 253--290, 
 II, 571--593, 
III, 809--829
\endref
\ref\key{\refLaxW}
\by \, \,   P.D. Lax and B. Wendroff
\paper Systems of conservation laws 
\jour  Comm. Pure Appl. Math. 
\vol 13
\yr 1960
\pages 217--237
\endref
\ref\key{\refMarcati}
\by \, \,   P.A. Marcati and R. Natalini
\paper Convergence of the pseudoviscosity approximation for
conservation laws
\jour Nonlin. Analysis T.M.A. \vol 23 \yr 1994  \pages 621--628
\endref
\ref\key{\refOleinik}
\by \, \,  O.A. Ole\v inik
\paper Uniqueness and Stability of the Generalized solution
of the Cauchy Problem for a Quasi-Linear Equation
\jour  Uspehi Mat. Nauk.
\vol14
\yr1959
\pages 165--170; {Amer. Math. Soc. Transl.} {33} (1963),  285--290
\endref
\ref\key{\refVonNeumann}
\by \, \,   J. Von Neumann and R.D. Richtmyer 
\paper A method for the
numerical calculation of hydrodynamical shocks
\jour J. Appl. Phys. 
\vol 21
\yr 1950
\pages 380--385
\endref
\ref\key{\refRaviart} 
\by \, \,  P.A. Raviart
\paper {Sur la r\'esolution num\'erique de
l'\'equation}
${\pa u\over\pa t}+u{\pa u\over\pa x}$$-\var{\alpha\over
\pa x}\left(\left|{\pa u\over\pa x}\right|{\pa u\over\pa x}\right)=0$
\jour J. Diff. Equa. \vol 8 \yr 1970 \pages 56-94
\endref
\ref\key{\refRichtmyer}
\by \, \,  R.D. Richtmyer and K.W. Morton
\book Difference methods for
initial-value problems\bookinfo 2nd ed.\publ Interscience Publ., J. Wiley and
Sons\publaddr New York\yr 1967
\endref
\ref\key{\refSchonbek}
\by \, \,   M.E. Schonbek
\paper Convergence of solutions to nonlinear dispersive equations
\jour Comm. Part. Diff. Equa.
\vol 7 
\yr 1982
\pages 959--1000
\endref
\ref\key{\refSzepessy} 
\by \, \, A. Szepessy
\paper An existence result for scalar conservation laws using measure-valued
solutions
\jour Comm. Part. Diff. Equa.
\vol 14
\yr 1989
\pages 1329--1350
\endref
\ref\key{\refTartar} 
\by \, \,  L. Tartar
\paper Compensated compactness and applications
to partial differential equations
\inbook Research Notes in Mathematics, Nonlinear
Analysis and Mechanics: Heriot-Watt Symposium, Vol. {\bf 4}, R. J.
Knops, New York, Pitman Press \yr 1979 
\pages 136--212
\endref
\ref\key{\refVenakidesOne}
\by \, \,  S. Venakides
\paper The zero dispersion limit of the periodic KdV equation
\jour  Amer. Math. Soc. Trans. 
\vol 301
\yr 1987
\pages  189--226
\endref
\ref\key{\refVenakidesTwo}
\by \, \,   S. Venakides
\paper  The Korteweg-de Vries equation with small dispersion:
higher order Lax-Levermore theory
\jour  Comm. Pure Appl. Math.
\vol 43
\yr 1990
\pages  335--361
\endref

\end{document}
